\magnification=1200

\centerline {\bf Embeddings of maximal tori in orthogonal groups}

\bigskip
\centerline
{Eva Bayer--Fluckiger}

\bigskip
Abstract : We  give necessary and sufficient conditions for an orthogonal group
defined over a global field of characteristic $\not = 2$ 
to contain a maximal torus of a given type.

\bigskip
{\bf Introduction}
\bigskip
Embeddings of maximal tori in orthogonal
groups have been studied in several papers, and occur in various arithmetic questions
(see for instance  [BCM], [F],  [G], [GR], [L], [PR]  and the references therein). The aim of
this paper is to give necessary and sufficient conditions for an orthogonal group defined
over a global field of characteristic $\not = 2$ 
to contain a maximal torus of a given type (see  Theorem  3.2.1). As we will see, this gives rise
to generalizations of 
some of the results of [F], [L] and [PR] (see Theorem 3.1.1 and Corollary 3.1.2).
\medskip
The case of tori of {\it type CM} (that is, tori associated to CM  \'etale algebras, see
1.2. and \S 4) is of special interest in some of the applications, and will be used here to
illustrate the results of the paper. 
The following is proved in \S 4~:
\medskip
\noindent
{\bf Theorem}. {\it Let $(E,\sigma)$ be a ${\bf Q}$--\'etale algebra  with involution of type CM of rank $2n$, and let $q$ be a quadratic
space over ${\bf Q}$ with ${\rm dim}(q) = {\rm rank}(E)$.  Then the orthogonal group
$O(q)$ contains a maximal torus of type $(E,\sigma)$ if and
only if ${\rm disc}(q) =  {\rm disc}(E) \in k^*/k^{*2}$,  the hyperbolicity condition holds {\rm (cf. 2.4)},
and the signature of $q$ is even.}

\medskip
In particular,  a torus of type CM can be embedded as a maximal torus of an orthogonal group if and only if such an embedding exists everywhere locally.

\bigskip
{\bf \S 1. Definitions, notation and basic facts}
\bigskip
Let $k$ be a field of characteristic $\not = 2$.
\medskip
{\bf 1.1.  Quadratic spaces}
\medskip
A {\it quadratic space} is a symmetric bilinear form of non--zero determinant
$q : V \times V \to k$, where $V$ is a finite dimensional $k$--vector space.
We denote by ${\rm dim}(q)$
its dimension (that is, the dimension of the underlying vector space $V$), and by
$O(q)$ its
orthogonal group. 
The {\it determinant} of  $q$ is
denoted by ${\rm det}(q)$; 
it is an element of $k^{\times}/k^{{\times}2}$. Let $m = {\rm dim}(q)$. Then the
{\it discriminant} of $q$ is by definition ${\rm disc}(q) = (-1)^{{m(m-1) \over 2}}{\rm det}(q)$. 
Let us denote by ${\rm Br}(k)$ the Brauer group of $k$,
considered as an additive abelian group, and let ${\rm Br_2}(k)$ be the subgroup
of elements of order $\le 2$ of ${\rm Br}(k)$. Any quadratic
space can be diagonalized, in other words there exist $a_1,\dots,a_m \in k^{\times}$ such
that $q \simeq <a_1, \dots ,a_m>$.  The {\it Hasse  invariant} of $q$ is 
by definition $\Sigma_{i < j} (a_i,a_j) \in {\rm Br_2}(k)$, where $(a_i,a_j)$ is the class
of the quaternion
algebra over $k$ determined by $a_i,a_j$, and is denoted by $w(q)$.  If $q$ and $q'$ are two quadratic spaces
over $k$, then we denote by $q \oplus q'$ their orthogonal sum. We have
$w(q \oplus q') = w(q) + w(q') + ({\rm det}(q),{\rm det}(q'))$ (see for instance [Sch, 2.12.6]).

\medskip
If $q : V \times V \to k$ is a quadratic space, 
let us denote by $\tau_q : {\rm End}(V) \to {\rm End}(V)$ the adjoint involution of $q$;
recall that we have $q(f(x),y) = q(x,\tau_q(f)(y))$ for all $f \in {\rm End}(V)$ and all $x,y \in V$. 

\bigskip
{\bf 1.2.  Maximal tori and \'etale algebras with involution}
\medskip
Recall that an \'etale algebra is a product of  separable field extensions of finite
degree of $k$. If $E$
is an \'etale algebra and $\sigma : E \to E$ is a $k$--linear involution, we denote by
$E^{\sigma}$ the subalgebra of $E$ fixed by $\sigma$.  The {\it unitary group}
$U(E,\sigma)$ is by definition the linear algebraic group over $k$ defined by
$U(E,\sigma)(A) = \{ x \in E \otimes_kA \ | \ x \sigma(x) = 1 \}$ for any commutative
$k$--algebra $A$. The following result is well--known (see for instance [BCM,  3.3], or [PR, 2.3]).
\medskip
\noindent
{\bf Proposition 1.2.1.}  {\it Let $q : V \times V \to k$ be a quadratic
space with ${\rm dim}(q) = 2n$.  Then we have

 \noindent {\rm (i)} 
Let $T \subset O(q)$ be a maximal $k$-torus. Then there is a unique \'etale
algebra $E \subset {\rm End}(V)$ stable by $\tau_q$ such that $T = U(E,\tau_q)$. 
Moreover, $E$ has rank $2n$ and $E^{\tau_q}$ has rank $n$.

\noindent {\rm (ii)}
Conversely, for any \'etale algebra $E \subset {\rm End}(V)$ stable under $\tau_q$ and satisfying
the rank  conditions above, the unitary group $U(E,\tau_q)$ is a maximal
$k$--torus of $O(q)$.  }

\medskip
If $q : V \times V \to k$ is a quadratic space and $E$ an \'etale algebra with involution $\sigma : E \to E$, 
we say that a maximal torus $T$ of $O(q)$ is of {\it type $(E,\sigma)$}  if  the conditions
of Proposition 1.2.1 hold for some \'etale algebra $E'  \subset {\rm End}(V)$ such
that the algebras with involution $(E,\sigma)$ and $(E',\tau_q|E')$ are isomorphic,  in particular $T \simeq U(E,\sigma)$.

\bigskip
{\bf 1.3.  Realizable pairs}
\medskip
If $(A_1,\tau_1)$ and $(A_2,\tau_2)$ are two $k$--algebras with involution. An {\it embedding}
of $(A_1,\tau_1)$ {\it  in}  $(A_2,\tau_2)$ is by definition
an injective homomorphism of algebras $A_1 \to A_2$ that commutes with the
involutions. 
\medskip
For any \'etale algebra with involution $(E,\sigma)$ and 
any $\alpha \in E^{\sigma}$, let $q_{\alpha} : E \times E \to k$
be the symmetric bilinear form $q_{\alpha}(x,y) = {\rm Tr}_{E/k}(\alpha x \sigma (y))$.
The following proposition is well--known
\bigskip
\noindent
{\bf Proposition 1.3.1.} {\it Let $(E,\sigma)$ be an \'etale algebra with involution of rank $2n$,
and assume that the rank of $E^{\sigma}$ is $n$. Let $q : V \times V \to k$  be a $2n$--dimensional
quadratic space. Then the following are equivalent :
\medskip
\noindent
{\rm (i)} 
The orthogonal group $O(q)$ contains a maximal torus of type $(E,\sigma)$.
\medskip
\noindent
{\rm (ii)}
The algebra with involution $(E,\sigma)$ can be embedded in the algebra with involution $({\rm End}(V),\tau_q)$.
\medskip
\noindent
{\rm (iii)}
There exists  $\alpha \in E^{\sigma}$
such that $q \simeq q_{\alpha}$.}
\medskip
\noindent
{\bf Proof.} The equivalence of {\rm (i)} and {\rm (ii)} follows from Proposition 1.2.1. For the
equivalence of {\rm (ii)} and {\rm (iii)}, see
for instance [PR, 7.1].
\bigskip
We say that the pair $(E,q)$ is {\it realizable} if the equivalent conditions of Proposition 1.2.1
hold. Recall that the discriminant of the \'etale algebra $E$ is by definition the determinant
of the quadratic space $E \times E \to k$ given by $(x,y) \mapsto {\rm Tr}_{E/k}(xy)$. 
It is denoted by ${\rm disc}(E)$. The following lemma is well--known, see for instance [F,  3.3.1] :

\bigskip
\noindent
{\bf Lemma 1.3.2.} {\it If $(E,q)$ is realizable, then ${\rm disc}(q) = {\rm disc}(E) \in k^{\times}/k^{\times 2}$.}
\medskip
\noindent
{\bf Proof.} As $q$ is realizable, we have $q = q_{\alpha}$ for some $\alpha  \in E^{\sigma}$.
Let $Q = q_1 : E \times E \to k$ and $Q' : E \times E \to k$ be the quadratic spaces defined by $Q(x,y) = {\rm Tr}_{E/k}(x \sigma (y))$ and $Q'(x,y) = {\rm Tr}_{E/k}(xy)$. We have
${\rm disc}(q) = {\rm N}_{E/k}(\alpha) {\rm disc}(Q)$. As $\alpha \in E^{\sigma}$, we
have ${\rm N}_{E/k}(\alpha) \in k^2$, hence ${\rm disc}(q) = {\rm disc}(Q)$. Writing
$E = E^{\sigma} (\sqrt \theta)$ for some 
$\theta \in E^{\sigma}$, a straightforward computation shows that 
${\rm det}(Q') = (-1)^n {\rm det}(Q)$.  As 
${\rm disc}(E) = {\rm det}(Q')$ and ${\rm disc}(Q) = (-1)^n {\rm det}(Q)$ by
definition, the result follows.

\bigskip
{\bf \S 2. Local conditions}
\medskip
Suppose that $k$ is a global field, and let us denote by $\Sigma_k$ the set of  places of $k$.
 We keep the notation of \S1.  Let $n \in {\bf N}$, and let $(E,\sigma)$ be an \'etale algebra with involution of rank $2n$.
Suppose  that  $E = K_1 \times \dots \times K_r$, where
$K_1, \dots, K_r$ are separable extensions of $k$, and that the $K_i$'s are all
stable by $\sigma$.  Let $I = \{1,\dots,r \}$, and for all
$i \in I$, let us denote by $F_i$ the fixed field of $\sigma$ in $K_i$. Suppose that $K_i$ is 
a quadratic extension of $F_i$ for all $i \in I$. Note that $E^{\sigma} = F_1 \times \dots \times F_r$,
and that ${\rm rank}(E^{\sigma}) = n$. Let $\Sigma_k^{\rm split}(E)$ be the set of $v \in \Sigma_k$ such that  all the
places of $E^{\sigma}$ above $v$ split in $E$.  
\bigskip 
We start by giving some local conditions for the embedding question of the previous
section.

\bigskip
{\bf 2.1. Split places}
\medskip
Recall that a quadratic
space $(V,q)$ is {\it hyperbolic} if there exists a subspace $W$ of $V$ such that
${\rm dim}(V) = 2 {\rm dim}(W)$, and $q(x,y) = 0$ for all $x,y \in W$. It is well--known
that a hyperbolic space is uniquely determined up to isomorphism by its dimension.
Let us denote by $h_{2n}$ the hyperbolic space of dimension $2n$.

\bigskip
\noindent
{\bf Lemma 2.1.1.} {\it Let $v \in \Sigma_k^{\rm split}(E)$,  and let $q$ be a $2n$--dimensional quadratic space over $k_v$.
Then  $(E,q)$ is realizable over $k_v$ if and only if  $q$ is hyperbolic.}
\medskip
\noindent
{\bf Proof.} As $v \in \Sigma_k^{\rm split}(E)$, over $k_v$ we have an isomorphism $E \simeq E_1 \times E_2$ where $E_1$
and $E_2$ are isomorphic \'etale $k_v$--algebras, and $\sigma(E_1) = E_2$. 
Let us show that  $q_{\alpha}$ is hyperbolic for any $\alpha \in (E^{\sigma}_v)^{\times}$.
Set $W = E_1 \times \{0 \}$. Then  $x \sigma (y) = 0$ for all $x,y \in W$, hence the restriction of $q_{\alpha}$ to $W$ is
identically zero. Since ${\rm dim}_k(W) = {1 \over2} {\rm dim}_k(E)$, this proves that $q_{\alpha}$ is hyperbolic, hence $q_{\alpha} \simeq h_{2n}$. Therefore $(E, h_{2n})$ is realizable over $k_v$. Conversely, if
$(E,q)$ is realizable over $k_v$, we have $q \simeq q_{\alpha}$ for some $\alpha \in (E^{\sigma}_v)^{\times}$, hence by the previous argument $q \simeq h_{2n}$.

\bigskip
{\bf 2.2. Non--split  places}
\bigskip
Recall that if $v \in \Sigma_k$ is a finite place or a real place, then ${\rm Br_2}(k_v)$ is a cyclic group of order 2. We will
identify it to $\{0,1 \}$. The following results will be used several times in the sequel. 
\bigskip
\noindent
{\bf Proposition 2.2.1.} {\it Let $v$ be a  place of $k$ such that $v \not \in \Sigma_k^{\rm split}(E)$.
Let $\epsilon \in \{ 0,1 \}$. Then there exists $\alpha \in (E_v^{\sigma})^{\times}$ 
such that $w(q_{\alpha}) = \epsilon$. }
\medskip
\noindent
{\bf Proof.} Recall that  $q_1 : E_v \times E_v  \to k_v$ is defined by $q_1(x,y) = {\rm Tr}_{E_v/k_v}(x \sigma(y))$.
If $w(q_1) = \epsilon$, we can take $\alpha = 1$. Suppose that $w(q) \not = \epsilon$. 
As $v \not \in \Sigma_k^{\rm split}(E)$, we have $E_v = E' \times K$, where $K$ is a field extension
of $k_v$ stable by $\sigma$. Set $F = K^{\sigma}$. Then $K$ is a quadratic extension of $F$.
Let $\beta \in F^{\times}$ such that $\beta \not \in {\rm N}_{K/F}(K^{\times})$. Let us denote
by $q_1'$ the restriction of $q_1$ to $K$. Then we have
$w(q_{\beta}) \not = w(q_1')$; this follows from
[M,  2.7] if $v$ is a finite place, and it
is clear if $v$ is an infinite place. Let $\alpha = (\beta,1) \in E_v^{\sigma}$.
Then $w(q_{\alpha}) \not = w(q_1)$, hence $w(q_{\alpha}) = \epsilon$.

 \medskip
\noindent
{\bf Lemma 2.2.2.} {\it Suppose that there exists a real place $u$ of $k$ such that 
we have $u \not \in \Sigma_k^{\rm split}(K_i)$
for all $i \in I$.
Then there exists a finite place $v$ of $k$  such that for all $i \in I$, we have $v \not \in \Sigma_k^{\rm split}(K_i)$.}

\medskip
\noindent
{\bf Proof.}  Let $L$ be a Galois extension of $k$ containing the fields $K_i$ for
all $i \in I$. Let $G = {\rm Gal}(L/k)$. Let us denote by $c$ the conjugacy
class of the complex conjugation in $G$ corresponding to an extension of
the place $u$ to $L$. By the Chebotarev density theorem,
there exists a finite place $v$ of $k$ such that the conjugacy class
of the Frobenius automorphism at $v$ is equal to $c$. Let $v$ be such a place.
Then all the places of $F_i$ above $v$ are inert in $K_i$. Therefore
we have $v \not \in \Sigma_k^{\rm split}(K_i)$ for all $i \in I$, and the
statement is proved.

\bigskip
{\bf 2.3. Real places}
\medskip
Let $v$ be a real place of $k$. It is well--known that any quadratic space $q$ over $k_v$ is isomorphic to 
$X_1^2 +  \dots + X_r^2 - X_{r+1}^2 - \dots - X_{r + s}^2$
for some non--negative integers  $r$ and $s$. These are uniquely determined by $q$, and we
have $r + s = {\rm dim}(q)$. The couple $(r,s)$ is called the {\it signature} of $q$ at $v$. 
We say that the signature of $q$ at $v$ is {\it even} if $r \equiv s \equiv 0 \  ({\rm mod} \ 2)$,
and we say that the signatures of $q$ are even if the signature of $q$ at $v$ is even
for all real places $v$ of $k$. 
\medskip
 We say that a place $w$ of $E^{\sigma}$ above $v$ is
{\it ramified} in $E$ if $w$ is a real place that extends to a complex place of $E$.
Let $\rho_v$ be the number of places of
$E^{\sigma}$ above $v$ which are not ramified  in $E$. The following lemma is
well--known

\bigskip
\noindent
{\bf Lemma 2.3.1.} {\it Let
$\alpha \in (E^{\sigma})^{\times}$. Then the signature of $q_{\alpha}$ is equal to
$(2 r_{\alpha}+ \rho_v, 2s_{\alpha} + \rho_v)$ where $r_{\alpha}$ is the number of
places of $E^{\sigma}$ above $v$  that ramify in  $E$ at which $\alpha$ is
positive,  and $s_{\alpha}$ is the number
places of  $E^{\sigma}$ that ramify in  $E$ at which $\alpha$ is
negative.}

\bigskip
\noindent
{\bf Proof.} This is immediate.

\bigskip
\noindent
{\bf Proposition 2.3.2.} 
 {\it Let
$q$ be a $2n$--dimensional quadratic space over $k_v$. 
Then $(E,q)$ is realizable if and only if the signature of $q$
is of the shape
$(2 r'+ \rho_v, 2s'+ \rho_v)$ for some non--negative integers
$r',s'$.}
\bigskip
\noindent
{\bf Proof.} If $(E,q)$ is realizable, then lemma 2.3.1.  shows that the signature
of $q$ has the required shape. Conversely, suppose that the signature of $q$
is equal to 
$(2 r'+ \rho_v, 2s'+ \rho_v)$ for some
$r',s' \in {\bf N}$. Let $\alpha \in (E^{\sigma})^{\times}$  be such that $\alpha$ is positive
at $r'$ places of $E^{\sigma}$ above $v$ and negative at $s'$ places. Then
by lemma 2.3.1, the signature of $q_{\alpha}$ is equal to $(2 r'+ \rho_v, 2s'+ \rho_v)$.
This implies that $q \simeq q_{\alpha}$, hence $(E,q)$ is realizable. 

\bigskip
{\bf 2.4. Combining local criteria}
\medskip If
$q$ is a $2n$--dimensional quadratic space over $k$,
we say that the {\it signature condition} holds for $E$ and $q$ if for every real
place $v$ of $k$, the signature of $q$
at $v$ is of the shape
$(2 r'+ \rho_v, 2s'+ \rho_v)$ for some non--negative integers
$r',s'$. For all $a \in {\rm Br}(k)$ and all $v \in \Sigma_k$,  let us  denote by $a_v$ the image of $a$ in
 ${\rm Br}(k_v)$.
Recall that $h_{2n}$ is the $2n$--dimensional hyperbolic space. We say that the {\it hyperbolicity condition} holds for $E$ and $q$ if for  all $v \in \Sigma_k^{\rm split}(E)$, we have 
$w(q)_v = w(h_{2n})_v$. 
\medskip
\noindent
{\bf Proposition 2.4.1.} {\it Let
$q$ be a $2n$--dimensional quadratic space over $k$. Then $(E,q)$ is
realizable over all the completions of $k$ if and only if ${\rm disc}(q) =  {\rm disc}(E) \in k^*/k^{*2}$,
and if the hyperbolicity condition and the signature condition hold for $q$ and $E$.}

\bigskip
\noindent
{\bf Proof.}  Suppose that ${\rm disc}(q) =  {\rm disc}(E) \in k^*/k^{*2}$,
and that  the hyperbolicity condition and the signature condition hold. Let us prove that
$(E,q)$ is realizable over $k_v$ for
all $v \in \Sigma_k$.  Suppose first that $v$ is an infinite place. If $v$ is complex, then there is nothing to prove. If $v$ is a real
 place, then by Proposition 2.3.2 the signature condition implies that $(E,q)$ is realizable
 over $k_v$. 
Suppose now that $v$ is a finite place. If $v \in \Sigma_k^{\rm split}(E)$,
 then the equality ${\rm disc}(q) =  {\rm disc}(E) \in k^*/k^{*2}$ and the
 hyperbolicity condition imply that the discriminants and the Hasse invariants of $q$ 
 and of $h_{2m}$ coincide over $k_v$. Therefore $q \simeq h_{2n}$ over $k_v$,
 and by Lemma 2.1.1  this implies that $(E,q)$ is realizable over $k_v$. 
 Suppose that $v \not  \in \Sigma_k^{\rm split}(E)$.
By Proposition  2.2.1, there exists $\alpha \in (E_v^{\sigma})^{\times}$ such that $w(q_{\alpha}) =
w(q)_v$. By Lemma 1.3.2, we have  ${\rm disc}(q_{\alpha}) = {\rm disc}(E)$. As by
hypothesis ${\rm disc}(q) =  {\rm disc}(E) \in k^*/k^{*2}$, the
discriminants of $q$ and $q_{\alpha}$ are equal in $k_v^{\times}/k^{\times 2}_v$. 
Therefore $q$ and $q_{\alpha}$ are isomorphic over $k_v$, and this implies that $(E,q)$ is realizable over $k_v$.  The converse follows immediately from Lemmas 1.3.2  and  2.1.1, and from
Proposition  2.3.2.

\bigskip
{\bf \S 3. Embedding criteria and Hasse principle}

\bigskip
We keep the notation of the previous sections. In particular, $k$ is a global field
of characteristic $\not = 2$, 
and 
$(E,\sigma)$ is \'etale algebra with involution of rank $2n$ such that  $E = K_1 \times \dots \times K_r$, where
$K_1, \dots, K_r$ are separable extensions of $k$,  the $K_i$'s are all
stable by $\sigma$, 
and  $F_i$ is  the fixed field of $\sigma$ in $K_i$ for all $i \in I = \{1,\dots,r \}$.

\bigskip  Recall that $\Sigma_k$ is the set of places of $k$, and that  $\Sigma^{\rm split}_k(K_i)$ is the set of $v \in \Sigma_k$  such that all
the places of $F_i$ above $v$ split in $K_i$. For all $i \not = j$, set $\Sigma_{i,j}=
\Sigma_k^{\rm split}(K_i) \cup \Sigma_k^{\rm split}(K_j)$. 
\bigskip
{\bf 3.1. Sufficient conditions and some notation}
\medskip
One of  the results of this section is  the following local--global principle

\medskip
\noindent
{\bf Theorem 3.1.1.} {\it  Suppose that there exists $i_0 \in I$ such that  for all $i \in I$, we
have $\Sigma_{i_0,i}  \not = \Sigma_k$. Let $q$ be a $2n$--dimensional quadratic space. Then a torus of type 
$(E,\sigma)$ can be embedded in the orthogonal group $O(q)$ 
if and only if such an embedding exists over all
the completions of $k$. }

\medskip Note that
this implies  [PR,  7.3] and [L, 2.20]. As we will see, Theorem 3.1.1 is a consequence of Theorem 3.2.1 below. We also get the following corollary, which provides an embedding criterion in
terms of invariants of the \'etale algebra and the quadratic space.

\medskip
\noindent
{\bf Corollary  3.1.2.} {\it 
Suppose that there exists $i_0 \in I$ such that  for all $i \in I$, we
have $\Sigma_{i_0,i}  \not = \Sigma_k$. Then 
$O(q)$ contains a maximal torus of type $(E,\sigma)$ if and
only if ${\rm disc}(q) =  {\rm disc}(E) \in k^*/k^{*2}$
and the signature and hyperbolicity conditions hold.}

\medskip
\noindent
{\bf Proof.}  
This follows from Proposition 2.4.1 and Theorem 3.1.1. 

\medskip
The following results will be needed in the proof of Theorem 3.1.1.

\bigskip
\noindent
{\bf Proposition 3.1.3.} {\it Suppose that $(E,q)$ is realizable over all the completions of $k$.
Then for all places $v$ of $k$ and $i \in I$, there exist quadratic spaces  $q^v_i$ over $k_v$ such that  
\smallskip
\noindent
{\rm (i)} for all $i \in I$ and every place  $v$ of $k$, the pair $(K^v_i,q^v_i)$ is realizable;
\smallskip
\noindent
{\rm (ii)} for every place  $v$ of $k$, we have $q \simeq  q^v_1 \oplus \dots \oplus  q^v_r$;
\smallskip
\noindent
{\rm (iii)} for all $i \in I$, we have $w(q^v_i) = 0$ for almost all $v \in \Sigma_k$. }

\bigskip
\noindent
Proposition 3.1.3 is an immediate consequence of Proposition  3.1.4  below, in which condition {\rm (iii)} 
is replaced by the more precise condition ${\rm (iii')}$. Let us start by introducing some
notation, that will be needed several times in the sequel. For all $i \in I$, let $n_i = [K_i:k]$, let 
 $d_i  = (-1)^{n_i} {\rm disc}(K_i)$, and set $D =  \Sigma_{i < j} (d_i,d_j) \in {\rm Br}_2(k)$. 
 Recall that for all $a \in {\rm Br}(k)$ and all $v \in \Sigma_k$, we denote by $a_v$ the image of $a$ in
 ${\rm Br}(k_v)$. Let $T$ be the set of places $v$ of $k$ such that $D_v \not = 0$, and let
$S$ be the set of places of $k$ at which the Hasse invariant of $q$ is
not equal to the Hasse invariant of the hyperbolic form of dimension equal to ${\rm dim}(q)$. 
Let $\Sigma_2$ be the set of dyadic places and $\Sigma_{\infty}$ the set of infinite places
of $k$, and
set $\Sigma = S \cup T \cup \Sigma_2 \cup \Sigma_{\infty}$. Note that $\Sigma$ is
a finite set. 
\medskip
\noindent
{\bf Proposition 3.1.4.} {\it Suppose that $(E,q)$ is realizable over all the completions of $k$.
Then for all places $v$ of $k$ and $i \in I$, there exist quadratic spaces  $q^v_i$ over $k_v$ such that  
\smallskip
\noindent
{\rm (i)} for all $i \in I$ and every place  $v$ of $k$, the pair $(K^v_i,q^v_i)$ is realizable;
\smallskip
\noindent
{\rm (ii)} for every place  $v$ of $k$, we have $q \simeq  q^v_1 \oplus \dots \oplus  q^v_r$;
\smallskip
\noindent
${\rm (iii')}$ for all $i \in I$, we have $w(q^v_i) = 0$ if  $v \not \in \Sigma$. }

\medskip
\noindent
{\bf Proof.} Let $v$ be a place of $k$. By hypothesis, $(E,q)$ is realizable over $k_v$. Hence
there exists $\alpha \in (E^{\sigma}_v)^{\times}$ such that $q \simeq q_{\alpha}$ over $k_v$, and
we have $\alpha = (\alpha_1,\dots,\alpha_r)$ with $\alpha_i \in (F_i^v)^{\times}$. Then the
quadratic spaces $q^v_i = q_ {\alpha_i}$ fulfill conditions {\rm (i)} and {\rm (ii)}. Let
us show that we can change the $q^v_i$ so that condition ${\rm (iii')}$  holds as well.
\medskip
Let $v \in \Sigma_k$ be such that $v \not \in \Sigma$, and suppose that there exists
$i \in I$ with $w(q^v_i) = 1$. Let us show that there exist quadratic spaces
$\tilde q^v_j$ for all $j \in I$ such that $w(\tilde q^v_j) = 0$ if $w( q^v_j) = 0$, and $w(\tilde q^v_i) = 0$.
As $v \not \in S \cup \Sigma_2$, we have $w(q)_v = 0$. Note that $w(q)_v = w(q_1^v) + \dots + w(q_r^v)  + D_v$,
and as $v \not \in T$, we have $D_v = 0$. Therefore there exists $m \in I$ with $m \not = i$
such that $w(q^v_m) = 1$. As $v$ is not dyadic, this implies that $q^v_i$ and $q^v_m$ are
not hyperbolic, hence by Lemma 2.1.1  we have $v \not \in \Sigma_{i,m}$.
As $v \not \in \Sigma_k^{\rm split}(K_i)$, by Proposition 2.2.1  there exists $\beta_i \in F_i^v$ such that
$w(q_{\beta_i}) = 0$. Similarly, as $v \not \in \Sigma_k^{\rm split}(K_m)$, there exists
$\beta_m \in F_m^v$ such that
$w(q_{\beta_m}) = 0$. 
Let $\tilde q^v_i = q_{\beta_i}$ and
$\tilde q^v_m = q_{ \beta_m}$, and set $\tilde q^v_j = q^v_j$ for $j \not = i,m$.  We have 
$w(\tilde q^v_j) = 0$ if $w( q^v_j) = 0$, and $w(\tilde q^v_i) = 0$.  By Lemma 1.3.2
we have ${\rm det}(\tilde q^v_j) = {\rm det}(q^v_j)$ for all $j \in I$. 
Moreover, as $w(\tilde q^v_i) = 0$ and $w(\tilde q^v_m) = 0$, we have $w(\tilde q^v_1 \oplus \dots \oplus \tilde q^v_r) =
w( q^v_1 \oplus \dots \oplus  q^v_r) $, implying that 
$\tilde q^v_1 \oplus \dots \oplus \tilde q^v_r  \simeq q^v_1 \oplus \dots \oplus  q^v_r$. Therefore condition {\rm (ii)} holds. The pairs ($K^v_j,\tilde q^v_j$) are realizable for all $j \in I$,
hence condition {\rm (i)} holds as well. Repeating this
procedure for all $i \in I$ with $w(q^v_i) = 1$ and for all $v \in \Sigma_k$ with $v \not \in \Sigma$ leads to
quadratic spaces over $k_v$ satisfying all three  conditions. This concludes the
proof of the proposition.

\bigskip
{\bf 3.2. A necessary and sufficient condition}
\medskip
In order to state a necessary and sufficient condition for the embedding problem
of tori in orthogonal groups (see Theorem 3.2.1 below), we need the following notation and definition
\medskip
\noindent
{\bf Notation.} 
Let ${\cal C}(E,q)$ be the set of collections $(q_i^v)$ of quadratic spaces over $k_v$ satisfying
conditions {\rm (i) - (iii)} of Proposition  3.1.3. 
For $C = (q_i^v) \in {\cal C}(E,q)$
and $i \in I$, 
set  $$S_i (C) = \{ v \in \Sigma_k'  \ | \ w(q^v_i) = 1 \}.$$ By condition {\rm (iii)} $S_i(C)$ is
a finite set, and we denote by $|S_i(C)|$ its cardinal.

\medskip
\noindent
{\bf Definition.} We say that $C = (q_i^v) \in {\cal C}(E,q)$ is {\it connected} if 
for all $i \in I$ such that  $|S_i(C)|$ is odd,
there exist $j \in I$ with $j \not = i$ such that  $|S_j(C)|$  is odd,
and  a chain $i = i_1, \dots, i_m = j$ of elements of $I$ with  $\Sigma_{i_t,i_{t+1}} \not = \Sigma_k$ for all $t = 1, \dots, m-1$.
 We say that  ${\cal C}(E,q)$ is
{\it connected} if it contains a connected element.

\medskip
\noindent
{\bf Theorem 3.2.1.} {\it Let $q$ be a $2n$--dimensional quadratic space. Then :

\noindent
{\rm (a)} The orthogonal group $O(q)$ 
contains a torus of type $(E,\sigma)$
over all completions of $k$ if and only if ${\cal C}(E,q)$ is not empty.

\noindent
{\rm (b)}  The orthogonal group $O(q)$ 
contains a torus of type $(E,\sigma)$
if and only if ${\cal C}(E,q)$ is connected.}

\medskip
\noindent
{\bf Proof.}  {\rm (a)}  With the terminology of 1.3, we have to show that
$(E,q)$ is realizable over all completions of $k$ if and only if ${\cal C}(E,q)$ not empty.
It is clear that if ${\cal C}(E,q)$ not empty, then $(E,q)$ is realizable over $k_v$ for 
all $v \in \Sigma_k$, and the
converse follows from Proposition 3.1.3. 
\medskip \noindent
{\rm (b)}  We have to prove that $(E,q)$ is realizable over $k$
if and only if ${\cal C}(E,q)$ is connected. If $(E,q)$ is realizable, then there
exist quadratic spaces $q_1,\dots,q_r$ over $k$ such that $q \simeq q_1 \oplus \dots \oplus q_r$
and that $(K_i,q_i)$ is realizable over $k$ for all $i \in I$. 
Set $q^v_i = q_i \otimes_k k_v$, and let $C = (q^v_i)$. Then $C  \in {\cal C}(E,q)$,
and $|S_i(C)|$ is even for all $i \in I$. Therefore $C$ is a connected element of ${\cal C}(E,q)$, hence ${\cal C}(E,q)$ is connected.
\medskip
Conversely, suppose that ${\cal C}(E,q)$ is connected,
and note that by part  {\rm (a)}  this implies that $(E,q)$ is realizable over all the completions of $k$.
Let us show that $(E,q)$ is realizable.
\medskip
\noindent {\bf Step 1.} If  $r = 1$, then $(E,q)$ is realizable. This can be deduced from [PR,  7.4] or [F, 1.1], but we give a (different) proof for the convenience of the reader. Let $v$ be a real place
of $k$ and let $(r_v,s_v)$ be the signature of $q$ at $v$. As $(E,q)$ is realizable over
$k_v$ by hypothesis, Proposition 2.3.2 implies that $(r_v,s_v) = (2r_v' + \rho_v, 2s'_v + \rho_v$)
for some $r'_v, s'_v \in {\bf N}$. Let $\alpha \in E^{\sigma}$ be such that $\alpha$ is positive at
exactly $r'_v$ real places of $E^{\sigma}$ that become complex in $E$. Then $\alpha$ is negative
at exactly real $s'_v$ places of $E^{\sigma}$ that become complex in $E$, hence by Lemma 2.3.1
the signature of $q_{\alpha}$ is $(r_v,s_v)$. Let $S_k$ be the set of places of $k$ at
which $q_{\alpha}$ and $q$ are not isomorphic. Note that $S_k$  consists of finite
places of $k$, and it is a finite set of even cardinality. If $v \in S_k$, then $v \not \in \Sigma^{\rm split}_k(E)$.
Indeed, both $(E,q)$ and $(E,q_{\alpha})$ are realizable over $k_v$ for all $v \in \Sigma_k$.
If $v \in \Sigma^{\rm split}_k(E)$, then by Lemma 2.1.1  this implies that $q$ and $q_{\alpha}$ are
both hyperbolic over $k_v$, hence they are isomorphic over $k_v$, and therefore
$v \not \in S_k$. For all $v \in S_k$, let us choose a place $w$ of $E^{\sigma}$ that does not
split in $E$ -- this is possible because $v \not \in \Sigma^{\rm split}_k(E)$. Let us denote
by $S_E$ the set of these places $w$. Then $S_E$ is in bijection with $S$, hence
it is also a finite set of even cardinality. Let us write $E = E^{\sigma}(\sqrt \theta)$ for
some  $\theta \in ({E^{\sigma}})^{\times}$, and let us choose $\beta \in (E^{\sigma})^{\times}$ such
that $(\beta,\theta)_w = -1$ if $w \in S_E$ and $(\beta,\theta)_w = 1$ otherwise.
This is possible as $S_E$ has even cardinality (see for instance
[O'M, 71.19], or [PR,  6.5]). Then by [M, 2.7],
the Hasse invariant of $q_{\alpha \beta}$ is equal to the Hasse
invariant of $q$. Since  these two quadratic spaces have equal dimension, determinant
and signatures, they are isomorphic by the Hasse--Minkowski theorem. Therefore
$(E,q)$ is realizable. 

\medskip
\noindent {\bf Step 2.} Let us show that ${\cal C}(E,q)$ contains  $C = (q_i^v) $ such that
\medskip
{\rm (iv)} {\it  $|S_i(C)| $ is even for all $i \in I$.}
\medskip
Let $C = (q_i^v) \in {\cal C}(E,q)$ be a connected element. Recall that by hypothesis $C$
satisfies conditions {\rm (i) - (iii)} of Proposition  3.1.3.
Suppose that for some $i \in I$, the integer $|S_i(C)|$ is odd. Since $C$ is connected, there exist
$j \in I$ with $j \not = i$ such that  $|S_j(C)|$  is odd,
and  a chain $i = i_1, \dots, i_m = j$ of elements of $I$ with  $\Sigma_{i_t,i_{t+1}} \not = \Sigma_k$ for all $t = 1, \dots, m-1$.  For all $t = 1, \dots, m-1$, let $v_t  \not  \in \Sigma_{i_t,i_{t+1}}$ be a finite place (note that this is possible by Lemma 2.2.2). Let $\alpha_1 \in (F^{v_1}_1)^{\times}$ be such that $q_i^{v_1} \simeq q_{\alpha_1}$ over $k_{v_1}$.
By Proposition 2.2.1, there exist $\alpha_t \in (F^{v_t}_t)^{\times}$
such that $w(q_{\alpha_t})  \not = w(q_{\alpha_{t+1}})$ for all $t = 1, \dots, m-1$. 
Set  $\tilde q^{v_t}_t = q_{\alpha_t}$ for all $t = 1, \dots, m-1$, and let
$\tilde q^u_s = q^u_s$ if $(u,s) \not = (v_t,t)$.  Set $\tilde C = (\tilde q_i^v)$. Then
$\tilde C  \in {\cal C}(E,q)$. We have $| S_i (\tilde C)|
 \equiv  0  \ \ ({\rm mod} \ 2)$, $| S_j (\tilde C)|
 \equiv  0  \ \ ({\rm mod} \ 2)$, and $| S_s (\tilde C)|
 \equiv  | S_s (C)|  \ \ ({\rm mod} \ 2)$ if $s \not = i,j$. 
Repeating this procedure we obtain
 a family of quadratic spaces safisfying conditions {\rm (i) - (iv)}.

\medskip
\noindent
{\bf Step 3.} End of proof. Let $C = (q_i^v) \in {\cal C}(E,q)$ satisfy conditions 
{\rm (i) - (iv)}; this is possible by Step 2.
For all $i \in I$, there exists a quadratic space $q_i$ over  $k$ such that
$q^v_i \simeq q_i$ over $k_v$ for all places $v$ of $k$.  This follows from [O'M, 72.1],
which applies because of conditions {\rm (iii)} and  {\rm (iv)}, and the fact that by condition  {\rm (i)} and
Lemma 1.3.2 we have ${\rm disc}(q^v_i) = d_i$
for all places $v$ of $k$. By condition {\rm (ii)}
we have $q \simeq q_1 \oplus \dots \oplus q_r$ over all the completions of $k$, 
hence by the Hasse--Minkowski theorem $q \simeq q_1 \oplus \dots \oplus q_r$
over $k$ as well. Note that by condition {\rm (i)}, the pair $(K_i,q_i)$ is realizable over all the completions of $k$.  
By Step 1, this implies that $(K_i,q_i)$ is realizable over $k$, hence $(E,q)$ is
realizable as well. This concludes the proof of the theorem.
\medskip
Note that the conditions (a) and (b) of Theorem 3.2.1 are not equivalent, in other
words the local--global principle does not hold in general : this follows from
the examples of Prasad and Rapinchuk, cf. [PR, 7.5].

\medskip
In order to deduce Theorem 3.1.1 from Theorem 3.2.1,  we need the following lemma
\medskip
\noindent
{\bf Lemma 3.2.2.} {\it Let $C = (q_i^v) \in {\cal C}(E,q)$. Then $\Sigma_{i \in I} | S_i (C)|
 \equiv  0  \ \ ({\rm mod} \ 2)$.}
 \medskip
 \noindent
 {\bf Proof.} For all $v \in \Sigma_k$, set $S_v(C) = \{ i \in I \ | \ w(q^v_i) = 1 \}$. 
We have $$\Sigma_{v \in \Sigma} |S_v(C)| = 
\Sigma_{i \in I} |S_i(C)|.$$ 
By  property {\rm (ii)}, we have
 $|S_v(C)| \equiv  w(q)_v  + D_v  \ \ ({\rm mod} \ 2)$ for all $v \in \Sigma_k$. 
 Therefore $\Sigma_{v \in \Sigma_k}  | S_v|(C)|
 \equiv  \Sigma_{v \in \Sigma_k}  w(q)_v +  \Sigma_{v \in \Sigma_k} D_v \ \ ({\rm mod} \ 2).$
 As $w(q)$ and $D$ are elements of ${\rm Br}_2(k)$, we have  $$\Sigma_{v \in \Sigma'_k}  w(q)_v 
  \equiv  0  \ \ ({\rm mod} \ 2), \ \ {\rm and} \ \ 
\Sigma_{v \in \Sigma'_k}  D_v   \equiv  0  \ \ ({\rm mod} \ 2).$$ This implies that
$\Sigma_{v \in \Sigma} | S_v(C)|
 \equiv  0  \ \ ({\rm mod} \ 2)$. As $\Sigma_{v \in \Sigma} |S_v(C)| = 
\Sigma_{i \in I} |S_i(C)|$, we also have $\Sigma_{i \in I} | S_i(C)|
 \equiv  0  \ \ ({\rm mod} \ 2)$. 

\bigskip 
\noindent
{\bf Proof of Theorem 3.1.1.} In order to apply Theorem 3.2.1, we have to show that ${\cal C}(E,q)$ is connected.
Let $C = (q_i^v) \in {\cal C}(E,q)$, and suppose that there exists $i \in I$ such that 
$|S_i(C)|$ is odd. By Lemma 3.2.2, we have $\Sigma_{i \in I} | S_i (C)|
 \equiv  0  \ \ ({\rm mod} \ 2)$. Therefore there exists $j \in I$ such that $j \not = i$, and that
$|S_j(C)|$ is odd. Since $\Sigma_{i_0,i}  \not = \Sigma_k$ and $\Sigma_{i_0,j}  \not = \Sigma_k$ by hypothesis, $C$ is
connected, and hence ${\cal C}(E,q)$ is connected. The result now follows from Theorem 3.2.1.

\medskip
Note that one can  give analogs of the results of  \S 3 in the
odd dimensional case. These can be easily deduced from the even dimensional case 
using the method of [PR, 7.2].

\bigskip
{\bf \S 4. An example - the case of CM \'etale algebras}
\bigskip
Recall that a number field is CM if it is a totally imaginary quadratic extension
of a totally real number field. Note that a number field is CM if and only if it has
exactly one complex conjugation (see for instance [Mi,  1.4.]). We say that $E$ is a CM \'etale algebra if it is a product of CM number fields, and the complex conjugation of $E$ is by
definition the product of the complex conjugations of its factors. 
\medskip
\noindent
{\bf Corollary 4.1.1.} {\it Suppose that 
$E$ is a CM \'etale algebra of rank $2n$, and that $\sigma : E \to E$ is the complex conjugation.
Let $q$ be a $2n$--dimensional quadratic
space over $k$. Then 
$O(q)$ contains a maximal torus of type $(E,\sigma)$ if and
only if  ${\rm disc}(q) =  {\rm disc}(E) \in k^*/k^{*2}$, the hyperbolicity condition holds
and signature of $q$ is even.}

\medskip
\noindent
{\bf Proof.} By Lemma 2.2.2, there exists $v  \in \Sigma_k$ such that for all $i \in I$, we have $v \not \in \Sigma_k^{\rm split}(K_i)$, Therefore for all $i,j \in I$ with $i \not = j$, we
have $\Sigma_{i,j}  \not = \Sigma_k$, and we can apply Corollary 3.1.2. 
As $E$ is CM and $\sigma$
is the complex conjugation,  we have
$\rho_v = 0$, hence the signature condition of Corollary 3.1.2  is equivalent to saying that the
signature of $q$ is even. 

\bigskip
{\bf Bibliography}
\bigskip \noindent
[BCM] R. Brusamarello, P. Chuard--Koulmann and J. Morales, Orthogonal groups containing
a given maximal torus, {\it J. Algebra} {\bf 266} (2003), 87--101. 
\medskip \noindent
[F] A. Fiori,  Special points on orthogonal symmetric spaces, {\it J. Algebra} {\bf 372} (2012),
397-419.
\medskip \noindent
[GR]  S. Garibaldi and A. Rapinchuk, Weakly commensurable S--arithmetic subgroups
in almost simple algebraic groups of types B and C, {\it Algebra and Number Theory},
to appear. 
\medskip \noindent
[G] P. Gille, Type des tores maximaux des groupes semi-simples, {\it J. Ramanujan
Math. Soc.}  {\bf 19} (2004), 213--230. 
\medskip \noindent
[L] T-Y. Lee, Embedding functors and their arithmetic properties, {\it Comment. Math. Helv}, to
appear. 
\medskip \noindent
[Mi] J. Milne, {\it Complex Multiplication}, http://www.jmilne.org/math/CourseNotes/cm.
\medskip \noindent
[M] J. Milnor, Isometries of inner product spaces, {\it  Invent. Math.} {\bf  8}  (1969), 83--97. 
\medskip \noindent
[O'M]  O.T. O'Meara, {\it Introduction to quadratic forms}, Reprint of the 1973 edition. Classics in Mathematics. Springer-Verlag, Berlin, 2000.
\medskip
\noindent
[PR] G. Prasad and A.S. Rapinchuk, Local--global principles for embedding of fields with
involution into simple algebras with involution, {\it Comment. Math. Helv.} {\bf 85} (2010),
583--645. 
\medskip \noindent
[Sch] W. Scharlau, {\it Quadratic and hermitian forms}, Grundlehren der Mathematischen Wissenschaften {\bf 270}, Springer-Verlag, Berlin, 1985.

\bigskip
\bigskip
Eva Bayer--Fluckiger

EPFL-FSB-MATHGEOM-CSAG

Station 8

1015 Lausanne, Switzerland

eva.bayer@epfl.ch

\bye